\documentclass{article}
\usepackage[latin1]{inputenc}
\usepackage[english]{babel}
\usepackage[T1]{fontenc}
\usepackage{amssymb}
\usepackage{amsmath,amsfonts,amscd,authblk}
\usepackage{mathtools}
\usepackage{graphicx}
\usepackage{indentfirst}
\usepackage{dsfont}
\usepackage{color}
\usepackage{enumitem}
\usepackage{float}
\usepackage{pict2e}
\usepackage{amsthm}
\usepackage{multirow}
\title{Efficiency of the averaged rank-based estimator for first order Sobol index inference}

\author[1]{Thierry Klein} 
\author[2]{Paul Rochet}

\affil[1]{Institut de Math\'ematiques de Toulouse; UMR5219. Universit\'e de Toulouse}
\affil[1,2]{ENAC - Ecole Nationale de l'Aviation Civile , Universit\'e de Toulouse, France.}

\newtheorem{theorem}{Theorem}[section]
\newtheorem{lemma}[theorem]{Lemma}
\newtheorem{prop}[theorem]{Proposition}

\newtheorem{coro}[theorem]{Corollary}

\newcommand{\var}{\operatorname{var}}
\newcommand{\cov}{\operatorname{cov}}

\newcommand{\opt}{{\operatorname{opt}}}
\newcommand{\rank}{{\operatorname{rank}}}

\date{}

\begin{document}

\maketitle

\abstract{Among the many estimators of first order Sobol indices that have been proposed in the literature, the so-called rank-based estimator is arguably the simplest to implement. This estimator can be viewed as the empirical auto-correlation of the response variable sample obtained upon re-ordering the data by increasing values of the inputs. This simple idea can be extended to higher lags of auto-correlation, thus providing several competing estimators of the same parameter. We show that these estimators can be combined in a simple manner to achieve the theoretical variance efficiency bound asymptotically.} \\

\textbf{Keywords:} Sensibility analysis, estimator averaging, asymptotic efficiency \\

\section{Introduction}

Sobol indices are by now a common tool for Global sensitivity methods which aim at detecting the most influential input variables/parameters in complex computer models. In this framework, the input variables are considered as random elements and   the relative  influence on the quantity of interest of each subset of its components is classically quantified by the Sobol indices, usually denoted by $S$  (see the book by Saltelli \cite{saltelli-sensitivity} for an overview on global sensitivity analysis). These indices, based on the Hoeffding's decomposition of the variance \cite{Hoeffding48}, were first introduced in  \cite{pearson1915partial} and later revisited in the framework of sensitivity analysis in \cite{sobol1993} (see also \cite{sobol2001global}). In a nutshell, a square integrable real-valued random variable $Y$, referred to as the \textit{output}, is entirely or partially explained by a collection $X$ of \textit{inputs} variables. The relative influence of $X$ on $Y$ is quantified by the Sobol index :
%$f$ be a $L^2$ mapping from $\mathbb{R}^p$ to $\mathbb{R}$, and $(X_i)_{1\leqslant i\leqslant p}$ be a collection of $p$ mutually independent random variables if $Y=f(X_1,\ldots,X_p)$, for any  subvector  $X$ of $(X_1,\ldots,X_p)$, the Sobol index with respect to $X$ is defined by
$$ S := \frac{\var(\mathbb E(Y | X) )}{\var(Y)} \in [0,1]. $$ 

In practice, an analytical expression of $S$ is rarely available making statistical inference on Sobol indices an important question. In the last decades, several approaches were developed in the literature, each one falling into one of the four following categories. 
\begin{itemize}
    \item Those  based on Monte Carlo,  quasi Monte Carlo  or nested Monte Carlo designs of experiments (see, e.g., \cite{Kucherenko2017different,Owen13,GODA201763}).
    \item  Those  based on spectral approaches (e.g.~Fourier Amplitude Sensitivity Test (FAST) \cite{cukier1978nonlinear}, Random Balance Design (RBD) \cite{tarantola2006random}, Effective Algorithm for computing global Sensitivity Indices (EASI).
    \cite{plischke2010effective} and polynomial chaos expansions \cite{sudret2008global} 
\item Those based on the so-called Pick freeze estimator in \cite{pickfreeze, janon2012asymptotic}.
\item Those based on a nearest neighbors approach \cite{devroye2003estimation,liitiainen2008nonparametric,
liitiainen2010residual,devroye2013strong,gyorfi2015asymptotic,devroye2018nearest, Chatterjee2019, GGKL20, broto2020variance}) or similar kernel-based methods \cite{zhu1996asymptotics,loubes2020rates}, studied in the particular case of first-order Sobol indices  \cite{da2009local,da2008efficient,plischke2020fighting,solis2021non,heredia2021nonparametric}
%,DVG} 
and in \cite{da2023new} for general Sobol indices.
\end{itemize}

Theoretical properties for the last two categories are well documented, especially in the case of first order Sobol indices. Consistency and asymptotic normality have been proved for kernel estimators \cite{da2008efficient, DVGKLP2023}, Pick Freeze \cite{janon2012asymptotic}, nearest neighbors estimators \cite{devroye2018nearest} as well as the rank based estimator \cite{GGKL20}. All these methods allow to estimate simultaneously all first-order Sobol indices from a single independent and identically distributed sample (two in the case of \cite{devroye2018nearest}), with the exception of the Pick freeze approach which requires a specific design of experiment associated to each input. \\

The kernel based approach developed in \cite{da2008efficient} is shown to be asymptotically optimal in quadratic mean, with its variance approaching the efficiency bound for a regular estimator of the conditional second order moment $\eta = \mathbb E \big( \mathbb E ( Y \, | \, X )^2 \big)$. However, the method is particularly tedious to implement and the estimator not easily tractable in practice. On the contrary, the rank-based approach developed by \cite{GGKL20} has by far the simplest implementation among all consistent methods but is sub-optimal in the sense that its variance does not reach the efficiency bound asymptotically. We show that the asymptotic variance of the rank estimator only differs from the efficiency bound by the additional term $\mathbb E \big( \var^2(Y | X) \big) $, which quantifies how far it is from optimality. \\

We introduce the family of lagged rank estimators $\widehat \eta^{(\ell)}, \ell \geq 1$ that generalizes the method of \cite{GGKL20}. We show that each lagged rank estimator $\widehat \eta^{(\ell)}$ performs similarly in quadratic mean as the original, under some control over the growth of the lag $\ell$ relative to the sample size $n$. By calculating the first order asymptotic expansion of the covariance matrix of a collection of lag estimators up to some maximal lag $k$, we derive an asymptotically optimal combination in the spirit of estimator averaging \cite{lavancier2016general}. More importantly, we show how the average estimator can be made to reach the efficiency bound of \cite{da2008efficient} by choosing $k$ growing sufficiently slowly to infinity relative to $n$.  \\

The article is organised as follows. We set the theoretical framework and the definition of the lagged rank estimators $\widehat \eta^{(\ell)}$ in Section \ref{sec:rank}. Their properties are investigated in Section \ref{sec:theory}, with a special focus on their joint second order moments and convergence in quadratic mean, paving the way to proving the efficiency of the averaging method. A numerical analysis to illustrate and validate the various results is presented in Section \ref{sec:num}. The proofs and technical lemmas are postponed to the Appendix.

\section{Rank estimators of Sobol indices}\label{sec:rank}

Let $(Y, X)$ be a couple of random variables with $Y$ real-valued and square-integrable. The Sobol index of $Y$ with respect to $X$, which measures the part of the variance of the output $Y$ that is ''explained'' by the input $X$, is given by
$$ S := \frac{\var(\mathbb E(Y | X) )}{\var(Y)} = \frac{\mathbb E \big( \mathbb E ( Y \, | \, X ) ^2 \big) - \mathbb E(Y)^2 }{\var(Y)}. $$ 
For inference purposes, because the expectation and variance of $Y$ do not depend on the input, the only real difficulty lies in estimating the second order conditional moment 
$$ \eta := \mathbb E \big( \mathbb E ( Y \, | \, X ) ^2 \big). $$
When $X$ is real-valued, in which case $S$ is generally referred to as a first-order Sobol index, a simple estimator of $\eta$ can be obtained from an iid sample $(Y_1, X_1), ..., (Y_n, X_n)$ following the method developed in \cite{GGKL20}. Let $(Y_{(i)}, X_{(i)})_{i=1,...,n}$ denote the data points sorted by increasing values of the $X_i$'s, i.e.~such that $ X_{(1)} \leq ... \leq X_{(n)}$, 
%and $\widehat m_Y , \widehat \sigma_Y^2$ the empirical mean and variance of $Y$. T
the rank estimator of $\eta$ is defined by
%$$ \widehat S = \frac{\frac 1 {n-1} \sum_{i = 1}^{n-1} Y_{(i)} Y_{(i+1)} - \widehat m_Y^2 }{\widehat \sigma_Y^2}.  $$
%Here, $\widehat m_Y $ and $\widehat \sigma_Y^2$ do not depend on the $X_i$'s and can be ignored in the estimation process, allowing to focus on the term
$$ \widehat \eta = \frac 1 {n-1} \sum_{i = 1}^{n-1} Y_{(i)} Y_{(i+1)} . $$
%which aims to approximate the second order conditional moment $ \eta := \mathbb E \big( \mathbb E ( Y \, | \, X ) ^2 \big)$. 
This estimator is known to be consistent and asymptotically Gaussian under mild conditions \cite{GGKL20}. A natural generalization of the idea consists in defining the lagged rank estimator associated to a lag $\ell \geq 1$ as
$$  \widehat \eta^{(\ell)} = \frac 1 {n-\ell} \sum_{i = 1}^{n-\ell} Y_{(i)} Y_{(i+\ell)} . $$
In order to investigate the properties of the lagged rank estimator $\widehat \eta^{(\ell)}$, let us introduce some technical assumptions related to the regularity of the relation between $Y$ and $X$. Let $\Phi, V$ be measurable functions such that $\Phi(X) = \mathbb E(Y | X)$ and $V(X) = \var(Y | X)$. We assume that $\Phi$ and $V$ are bounded
\begin{equation}\label{H1}\tag{H1} \forall x \in \mathbb R \, , \,  | \Phi(x) | \leq M_\Phi  \, \text{ and }  | V(x)  | \leq M_V
\end{equation}
and Lipschitz
\begin{equation}\label{H2}\tag{H2} \forall x, x' \in \mathbb R \, , \,  | \Phi(x) - \Phi(x')  | \leq L_\Phi |x - x' |  \, \text{ and }  | V(x) - V(x')  | \leq L_V |x - x' | 
\end{equation}
for some positive constants $M_\Phi, M_V, L_\Phi, L_V$. Remark that under these assumptions, $\Phi^2$ is also bounded 
%by $M_\Phi^2$ 
and Lipschitz, 
%with Lipschitz constant $2 M_\Phi L_\Phi $
which will be useful in later proofs. These conditions are quite mild compared to the usual assumptions for first order Sobol index inference. For instance, it is extremely common in the literature to assume that the inputs are uniformly distributed on $[0,1]$ or have compact support. In this case, it is typically sufficient to assume that the conditional expectation and variance are continuously differentiable for \eqref{H1} and \eqref{H2} to hold. \\

%Let $\mathcal F_n = (X_1,...,X_n)$ and denote by $\Delta_{n}$ the sample range $\Delta_{n} = X_n - X_1$. 
%$$ \Delta_{k, n} = \Delta_{k, n}(\mathcal F_n) : = \displaystyle \sup_{i=1,..., n-2k} | X_i - X_{i+2 k} | . $$
%$$ \Delta_{k, n} = \Delta_{k, n}(\mathcal F_n) : = \frac {1} {n-k} \sum_{i=1}^{n-k}| X_{i+k} - X_i |. $$
In the sequel, we shall denote by $k = k(n)$ the total number of lags considered for our purposes. This number is allowed to increase with $n$ but must somehow be constrained by the distribution of the inputs, in particular by their range 
$$ \Delta_n : = X_{(n)} - X_{(1)} . $$
Typically, we want to be able to consider as many lags $\ell$ as possible provided that the average distance between two data points $X_{(i)}$ and $X_{(i + \ell)}$ is sufficiently small (roughly speaking, we want $X_{(i)}$ to be close enough to $X_{(i+ \ell)}$ so that $Y_{(i)}$ and $Y_{(i+\ell)}$ are almost identically distributed conditionally to $X_{(i)}$). This way, $\widehat \eta^{(\ell)}$ should provide an accurate depiction of the second conditional moment of $Y$. By a telescoping argument, the average distance can be bounded by
\begin{equation}\label{eq:teles}  \frac 1 {n-\ell} \sum_{i=1}^{n-\ell} \big( X_{(i + \ell)} - X_{(i)} \big) \leq \frac{\ell }{n-\ell} \Delta_n \leq \frac{k }{n-k} \Delta_n .
\end{equation}
%\frac{\ell }{n-\ell} \big( X_{(n)} - X_{(1)} \big) 
%\leq \frac{k}{n-k}| X_n - X_1 | 
%=: \Delta_{\ell,n} $$
Hence, we require this term to vanish fast enough as $n \to \infty$, via the following simple assumption
%\begin{equation}\label{H3}\tag{H3} k = k_n = o \big(\sqrt n \big)  \, \text{ and }  \, \mathbb E \big(\Delta_{k,n}^2 \big) = o \Big( \frac 1 n \Big).
%\end{equation}
\begin{equation}\label{H3}\tag{H3} \mathbb E \big(k^2 \Delta_{n}^2 \big) = o ( n).
\end{equation}

This condition can be understood as both a regularity assumption on the tail of the input's distribution and a restriction on the maximal number $k = k(n)$ of lags considered. It is nonetheless quite mild and can always be met unless the distribution of the inputs is heavy tailed, leading to extreme behaviors of the inputs' range $\Delta_n $. The minimal requirement, corresponding to the situation where the  distribution of the $X_i$'s has compact support (excluding the trivial case $\Delta_n \overset{a.s.}{=} 0$), is to take $k = o(\sqrt n)$. 
%For a Gaussian distribution, the sample range $|X_n -X_1 |$ grows as $ \sqrt{\log n}$, requiring to choose $k = o(\sqrt{n/\log n})$, a value far from restrictive in practice. 
If the distribution of the inputs decays exponentially fast, the asymptotic behavior $ \Delta_n = O_P(\log n)$ imposes the slightly stronger condition $k = o \big(\sqrt n / \log n \big)$, far from prohibitive in practice. Finally, remark that we do not rule out data-driven values of $k$, such as e.g.~$k \sim n^{1/3}/\Delta_n$ which automatically satisfies \eqref{H3} regardless of the distribution of the inputs. Nevertheless, the cautious and simple $k = \lfloor n^{1/3} \rfloor$, which we use in all numerical applications, fulfills all theoretical requirements while providing a good rule of thumb for practical purposes, as discussed in Section \ref{sec:num}. 

\section{Theoretical results}\label{sec:theory}

We are now in position to investigate some properties of the lagged rank estimators. Because we are ultimately interested in their convergence in quadratic mean, we focus on controlling the bias and variance, both for a finite sample size $n$ and asymptotically as $n$ grows to infinity. Only the main results are presented in this section, the detailed proofs and technical steps can be found in the Appendix.

\begin{prop}\label{prop:exp} Under \eqref{H1}, \eqref{H2}, we have for all $\ell = 1,...,k$, 
$$ \big| \mathbb E \big( \widehat \eta^{(\ell)} \big) - \eta  \big| \leq \frac{\ell}{n-\ell} \Big( L_\Phi M_\Phi  + 2 M_\Phi^2 \, \mathbb E \big( \Delta_n \big)  \Big).  $$
In particular, if \eqref{H3} is also met, then $ \mathbb E \big( \widehat \eta^{(\ell)} \big) = \eta +  o \big( n ^{-1/2} \big) $.
\end{prop}

This result, which is a direct consequence of Lemma \ref{lem:espX} in the Appendix, illustrates how the bias of $\widehat \eta^{(\ell)}$ may strongly depend on the lag $\ell$. We observe this phenomenon in some examples of the numerical analysis in Section \ref{sec:num} where the bias term is shown to highly vary in function of the lag, especially for smaller sample sizes $n$. Nevertheless, the variance becomes the dominating term asymptotically, as shown in the next proposition. \\

\begin{prop}\label{prop:var} Under \eqref{H1}, \eqref{H2} and \eqref{H3}, we have for all $\ell = 1,...,k$, 
$$ n \var \big( \widehat \eta^{(\ell)} \big) = 4 \, \mathbb E \big( \Phi^2(X) V(X) \big) + \mathbb E \big(V^2(X) \big) + \var \big(\Phi^2(X) \big)  + o (1 ) .  $$
%\mathbb P_X \big( g_2^2 \big) + 2 \mathbb P_X \big( g_2 \Phi^2 \big)  - 2 \mathbb P_X \big(\Phi^4 \big) - \mathbb  P_X (\Phi^2)^2  + o (1 ) .  $$
\end{prop}

Let us compare the limit variance to that of other existing estimators of single input Sobol indices. The main term (up to the convergence rate of $1/n$), given by
$$ \sigma^2_{\operatorname{rank}} = 4 \, \mathbb E \big( \Phi^2(X) V(X) \big) + \mathbb E \big(V^2(X) \big) + \var \big(\Phi^2(X) \big)  $$
falls short to the theoretical optimal value
$$  \sigma^2_{\operatorname{opt}} = 4 \, \mathbb E \big( \Phi^2(X) V(X) \big) + \var \big(\Phi^2(X) \big) $$
shown in \cite{da2008efficient} to be the asymptotic lower bound for the variance of an estimator of $\eta$. In the same paper, the authors propose a method that achieves the theoretical lower bound for the asymptotic variance, but relies on a preliminary non-parametric estimation of the joint density of $(X,Y)$  along with various tuning parameters, making its construction somewhat tedious. Note that the rank estimator $\widehat \eta^{(\ell)}$ is asymptotically optimal if, and only if, $V(X) \overset{a.s.}{=} 0$
%, that is, if $ Y \overset{a.s.}{=} \Phi(X)$ 
in which case the Sobol index is equal to one.  \\

For the sake of comparison, the alternative estimator of $\eta$ proposed in \cite{devroye2018nearest} and based on a nearest neighbors estimation of the conditional expectation, achieves an asymptotic theoretical variance of
$$ \sigma^2_{\operatorname{nn}} =   5 \, \mathbb E \big( \Phi^2(X) V(X) \big) + 2 \mathbb E \big(V^2(X) \big) + 2 \var \big(\Phi^2(X) \big) .       $$
While the three variances are always comparable, 
$$ \sigma^2_{\operatorname{opt}} \leq \sigma^2_{\operatorname{rank}} \leq \sigma^2_{\operatorname{nn}}, $$
an important advantage of the nearest neighbors approach over the rank method is that it can handle the estimation of multiple inputs Sobol indices, a problem that notoriously suffers from the curse of dimensionality. More recently, a kernel approach inspired from \cite{GGKL20} was proposed in \cite{DVGKLP2023} with an asymptotic theoretical variance of
$$ \sigma^2_{\operatorname{ker}} =   4 \, \mathbb E \big( \Phi^2(X) V(X) \big) + 4  \var \big(\Phi^2(X) \big) .       $$
This variance is, of course, higher than the theoretical lower bound $\sigma^2_{\operatorname{opt}}$ but is not comparable to the other two. \\
%Nevertheless, as it can be easily seen depending on the example that  $\sigma^2_{\operatorname{ker}}$ and $\sigma^2_{\operatorname{rank}}$ (resp. $\sigma^2_{\operatorname{ker}}$ and $\sigma^2_{\operatorname{D}}$) are not comparable. 

From an implementation point of view, the rank estimator $\widehat \eta^{(\ell)}$ is by far the easiest to construct, with the ordering of the inputs as its main computational hurdle. Besides its simplicity, a notable advantage of the method is to provide a new estimator for each lag $\ell$, with similar properties asymptotically. This feature can be exploited by combining an appropriate number of rank estimators obtained with different lags, in order to improve the estimation. The next result shows how the rank estimators $\widehat \eta^{(1)}, ..., \widehat \eta^{(k)}$ form a collection of competing estimators with symmetric behaviors asymptotically.

\begin{prop}\label{prop:cov} Under \eqref{H1}, \eqref{H2} and \eqref{H3}, we have for $1 \leq \ell < m \leq k$, 
$$ n  \cov \big( \widehat \eta^{(\ell)}, \widehat \eta^{(m)} \big) = 
4 \, \mathbb E \big(\Phi^2(X) V(X) \big) + \var \big(\Phi^2(X) \big) 
+ o (1).  $$
\end{prop}

\noindent For a fixed $k$, Propositions \ref{prop:var} and \ref{prop:cov} give the following first order term in the asymptotic expansion of the covariance matrix $\Sigma := \big( \Sigma_{\ell m} \big)_{\ell, m =1,...,k} $ of $\widehat \eta^{(1)}, ..., \widehat \eta^{(k)}$ :
\begin{equation}\label{eq:sigma_asympt}  \Sigma_{\ell m} = \lim_{n \to \infty} n \, \cov \big( \widehat \eta^{(\ell)}, \widehat \eta^{(m)} \big) = \left\{ \begin{array}{cl} \sigma^2_{\operatorname{opt}}  + \mathbb E \big(V^2(X) \big) & \text{ if } \ell = m \\ \sigma^2_{\operatorname{opt}} & \text{ if } \ell \neq m \end{array} \right. % \, \, , \, \, \ell, m = 1,...,k. 
\end{equation}
Remark that $\Sigma$ is of full rank provided that $V(X)$ is not almost surely zero, which indicates that the $\widehat \eta^{(\ell)}$'s are linearly independent asymptotically. Therefore, it is possible to reduce the asymptotic variance of an estimator of $\eta$ by considering a linear combination 
$$ \widehat \eta^{(k)}_{\operatorname{av}} = \sum_{\ell = 1}^k \lambda_\ell \, \widehat \eta^{(\ell)}, $$
where the weights $\lambda_\ell, \ell = 1,...,k$ are constrained to sum up to one. This heuristics is investigated in \cite{lavancier2016general} to determine the weights minimizing the asymptotic variance as a function of $\Sigma$. Although $\Sigma$ is unknown in practice, the symmetrical form of $\Sigma$ in this case, having the same diagonal values as well as off-diagonal values, suffices to deduce that the solution corresponds to the equal weights $\lambda_\ell = 1/k$. This simple way of combining the rank estimators actually achieves the theoretical efficiency bound of \cite{da2008efficient} under mild assumptions, as shown in the next theorem. \\

\begin{theorem}\label{th:av} If $k = k(n)$ tends to infinity as $n \to \infty$ and the conditions \eqref{H1}, \eqref{H2} and \eqref{H3} are met, the average estimator $\widehat \eta^{(k)}_{\operatorname{av}} $ obtained with equal weights $\lambda_\ell = 1/k$ satisfies
%$\sqrt n \, \big( \eta - \mathbb E ( \widehat \eta^{(k)}_{\operatorname{av}}) \big) = o (1) $ and
%$$ n \var \big( \widehat \eta^{(k)}_{\operatorname{av}} \big) = 4  \mathbb E \big( \Phi^2(X) \sigma^2(X) \big) + \var \big(\Phi^2 (X) \big) + \frac 1 k \mathbb E \big( \sigma^4(X) \big) +  o(1).  $$
%In particular, if $k = k_n$ tends to infinity sufficiently slowly so as to still verify the assumptions of Proposition \ref{prop:var}, then $\widehat \eta^{(k)}_{\operatorname{av}}$ is asymptotically efficient with asymptotic variance
$$ \lim_{n \to \infty}  n \var \big( \widehat \eta^{(k)}_ {\operatorname{av}} \big) = 4 \, \mathbb E \big(\Phi^2(X) V(X) \big) + \var \big(\Phi^2(X) \big) = \sigma^2_{\operatorname{opt}}.   $$
\end{theorem}

The fact that the averaged rank estimator achieves the variance efficiency bound $\sigma^2_{\operatorname{opt}}$ as $n \to \infty$ is certainly encouraging, although the result concerns the actual mean square error (MSE) of the estimator and not the variance of the Gaussian limit for a regular estimator, as introduced in \cite{da2008efficient}. While the regularity and asymptotic normality of $\widehat \eta^{(k)}_{\operatorname{av}} $ are to be expected under the appropriate assumptions, it has not been investigated in this paper as it deviates from the original objective of variance reduction.

\section{Numerical analysis}\label{sec:num}
We investigate the performances of the rank estimators and their averages in different models of the form
$$ Y = \Phi(X) + \sqrt{V(X)} \, \epsilon  $$
where $\epsilon$ is a standard Gaussian random variable independent from $X$. Each model is simulated $N = 10000$ times to give a faithful representation of the distributions of the different estimators. We show the boxplots of the rank estimators obtained for all lags from $\ell = 1$ to $\ell = 50$ for four samples sizes from $n=100$ to $n=2000$, which we compare to the boxplots of the averages obtained for $k=5$ to $k=50$ with $5$ estimators added at each step. \\

Due to the similarities in the interpretations of the results produced from various models, we choose to discuss only two values of the conditional expectation function, namely $\Phi(X) = \sin(5X) $ and $\Phi(X) = X^2 - 3X$. For the conditional variance, all examples are generated with $V(X) = 4X^2$ as other values of $V$ hardly had any noticeable impact on the results. For the distribution of the inputs $X_i$, we considered the uniform distribution on $[0,1]$ and the standard exponential distribution. 

\begin{figure}[H]
\centering
\includegraphics[width = 0.95\textwidth]{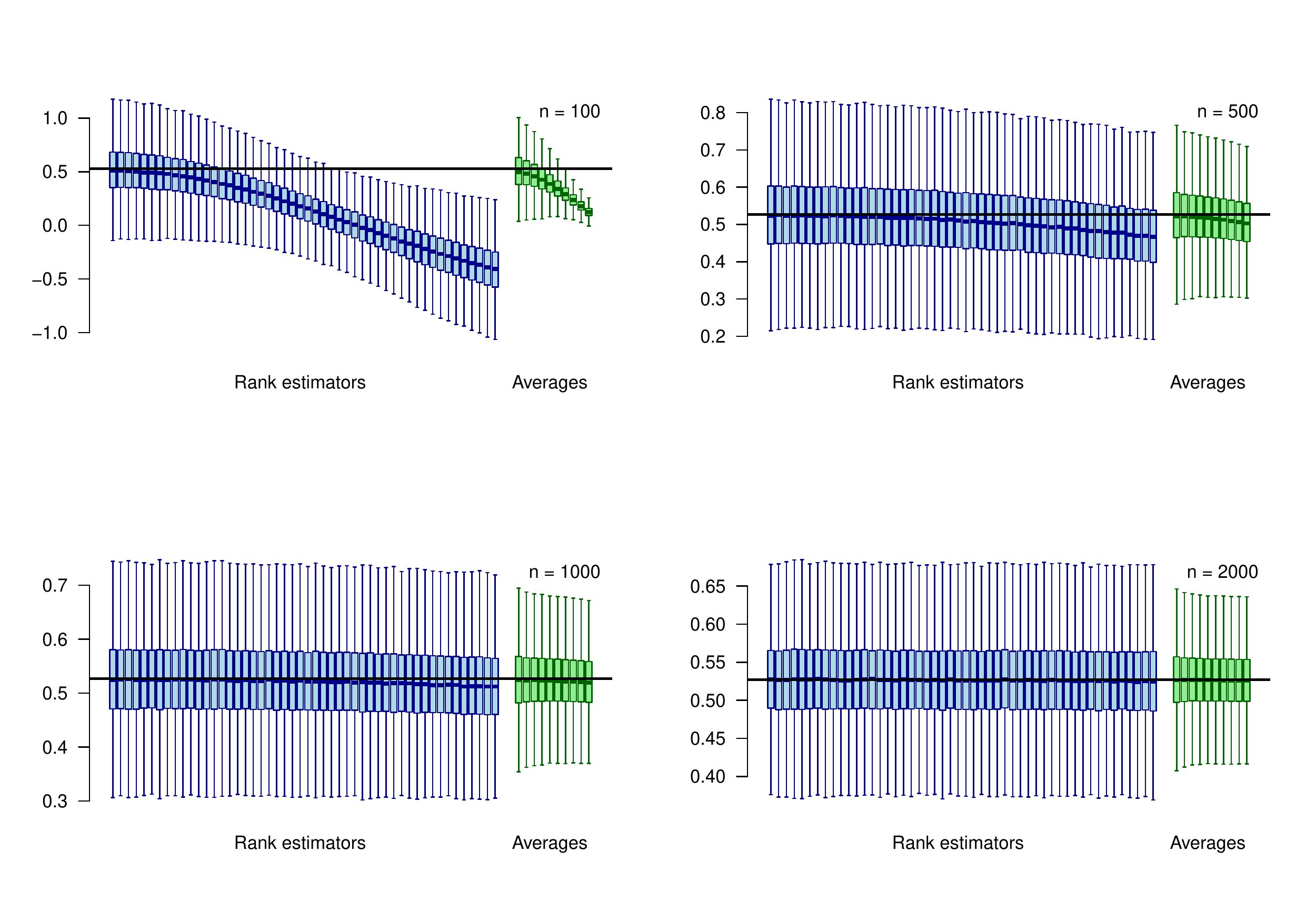}
\caption{\footnotesize The model $Y = \sin(5X) + 2X \epsilon$ where $X$ is uniformly distributed on $[0,1]$  and $\epsilon$ is a standard Gaussian random variable independent from $X$. In blue, the boxplots of the rank estimators of $\eta$ from lag $\ell = 1$ to $\ell = 50$. In green, the average estimators obtained with $k=5$ to $k=50$ by steps of $5$. The sample sizes vary from $n=100$ (top-left) to $n = 2000$ (bottom-right). The boxplots are constructed by Monte-Carlo using $N=10^4$ repetitions.}
\label{fig:sin_unif}
\end{figure}

In Figure~\ref{fig:sin_unif}, we observe that the bias of the rank estimators is important and varies strongly with the lag for the smaller sample sizes $n$, but does vanish asymptotically as predicted by the theory. The averaging procedure appears to improve significantly the performances of the rank estimators, as can be expected in this model with a maximal theoretical improvement of around $(\sigma^2_\rank - \sigma^2_\opt)/\sigma^2_\rank \approx 49 \%$. The positive effect of the averaging is mostly visible on the variance (smaller inter-quartile intervals) but can not compensate for the biases, all of the same sign.  

\begin{figure}[H]
\centering
\includegraphics[width = 0.95\textwidth]{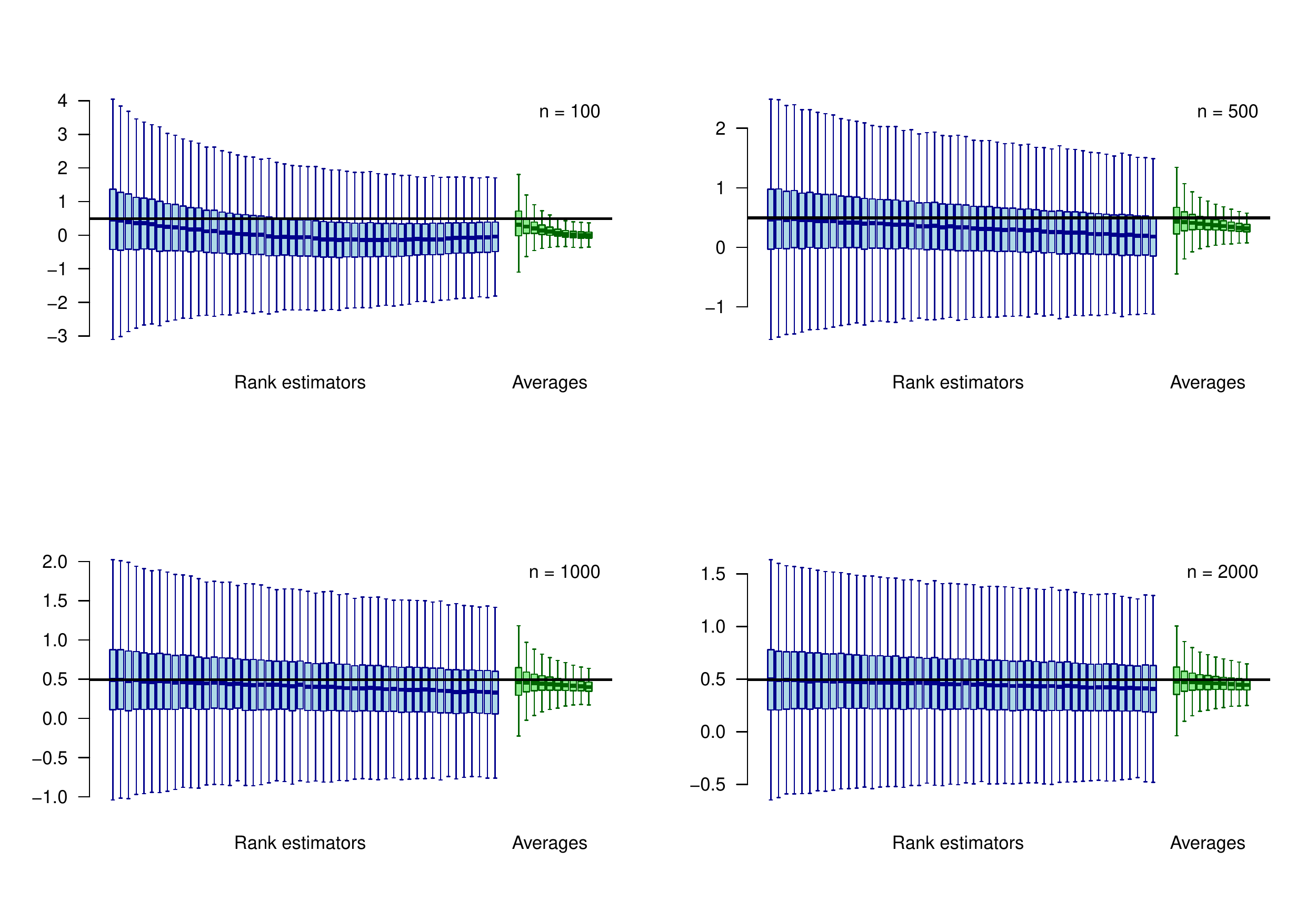}
\caption{\footnotesize The model $Y = \sin(5X) + 2X \epsilon$ where $X$ has standard exponential distribution and $\epsilon$ is a standard Gaussian random variable independent from $X$. In blue, the boxplots of the rank estimators of $\eta$ from lag $\ell = 1$ to $\ell = 50$. In green, the average estimators obtained with $k=5$ to $k=50$ by steps of $5$. The sample sizes vary from $n=100$ (top-left) to $n = 2000$ (bottom-right). The boxplots are constructed by Monte-Carlo using $N=10^4$ repetitions.}
\label{fig:sin_exp}
\end{figure}

A similar behaviour can be observed in Figure~\ref{fig:sin_exp} despite the distribution of the input $X$ not having compact support. The maximal theoretical improvement from the averaging procedure is even higher in this case, being around $(\sigma^2_\rank - \sigma^2_\opt)/\sigma^2_\rank \approx 96 \%$. \\

\begin{figure}[H]
\centering
\includegraphics[width = 1\textwidth]{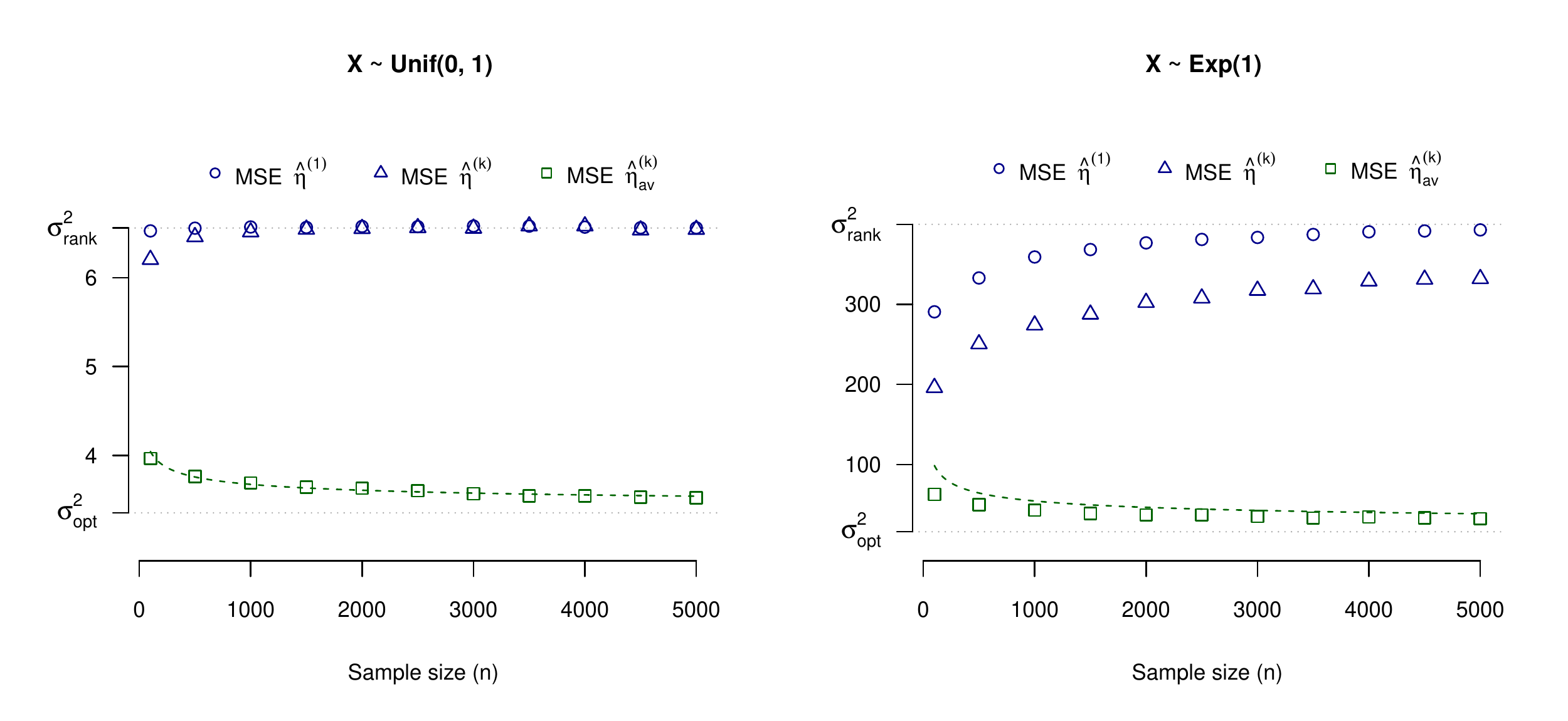}
\caption{\footnotesize Evolution of the MSE's of the estimators $\widehat \eta^{(1)}$, $\widehat \eta^{(k)}$ and the average $\widehat \eta_{\operatorname{av}}^{(k)}$ with $k = \lfloor n^{1/3} \rfloor$, as a function of $n$ in the model $Y = \sin(5X) + 2X \epsilon$, under uniform (left) and exponential (right) distribution of the inputs. The MSEs are multiplied by the sample size $n$ to highlight the convergence as $n \to \infty$. The values are obtained from Monte-Carlo estimations with $N = 10^5$ repetitions. The green dashed line represents the MSE first order asymptotic expansion for the averaged estimator derived from Equation~\eqref{eq:asympt_k}.}
\label{fig:vars}
\end{figure}

The convergence in quadratic mean of the rank and averaged estimators are sensible to the regularity conditions of the model, as can be seen in Figure~ \ref{fig:vars}. In the model $Y = \sin(5X) + 2X \epsilon$, with uniformly distributed inputs, where the regularity conditions \eqref{H1} and \eqref{H2} are satisfied, the MSEs of the various estimators do appear to behave accordingly to the theory in function of the sample size, rapidly reaching the asymptotic regime. The numerical results are not as convincing in the same model with exponentially distributed inputs, where the various estimators are slower to reach their asymptotic regime. This is especially true for the lagged rank estimator $\widehat \eta^{(k)}$ with $k$ growing to infinity, although it surprisingly performs better than expected by the theory. Remark that in this case, none of the conditions \eqref{H1} and \eqref{H2} hold for the conditional variance $V$, which is neither bounded nor Lipschitz on the support of the inputs distribution. Nevertheless, the evolution of the MSE of the averaged estimator seems to validate in both cases the theoretical first order expansion
$$  n \, \var \big( \widehat \eta^{(k)}_{\operatorname{av}} \big) \approx \sigma^2_{\operatorname{opt}} + \frac 1 k  \mathbb E \big(V^2(X) \big) $$
derived from Equation~\eqref{eq:asympt_k} in the proof of Theorem \ref{th:av}. In all these scenarios, the squared bias account for less than $1\%$ of the MSE, making it indistinguishable from the variance in the graphical representations. \\

\begin{figure}[H]
\centering
\includegraphics[width = 0.95\textwidth]{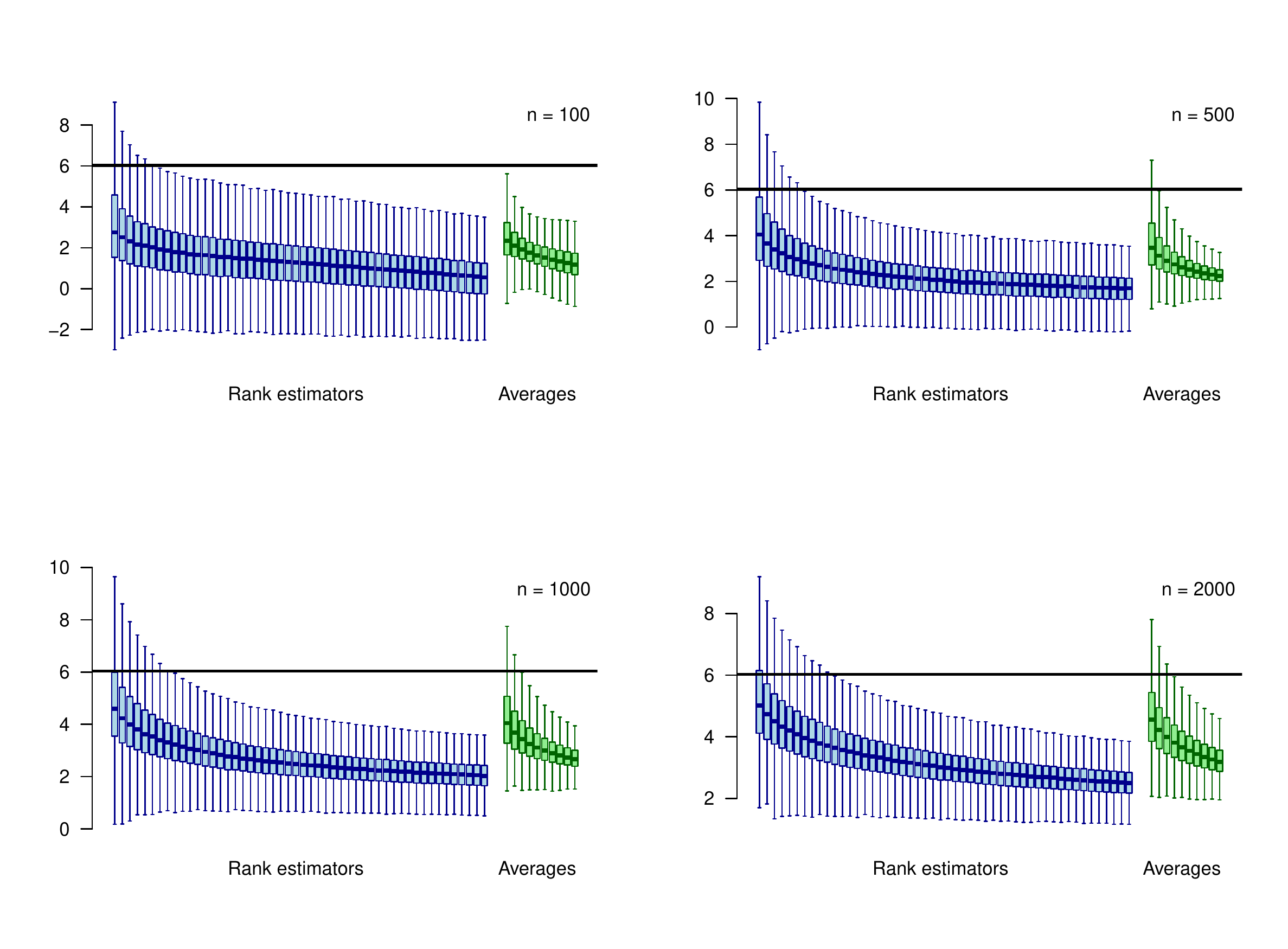}
\caption{\footnotesize The model $Y = X^2 - 3X + 2X \epsilon$ where $X$ has exponential distribution and $\epsilon$ is a standard Gaussian random variable independent from $X$. In blue, the boxplots of the rank estimators of $\eta$ from lag $\ell = 1$ to $\ell = 50$. In green, the average estimators obtained with $k=5$ to $k=50$ by steps of $5$. The sample sizes vary from $n=100$ (top-left) to $n = 2000$ (bottom-right). The boxplots are constructed by Monte-Carlo using $N=10^4$ repetitions.}
\label{fig:x2_exp}
\end{figure}

Figure~\ref{fig:x2_exp} illustrates how things can fall apart when the regularity conditions in \eqref{H1} and \eqref{H2} are not met for the conditional expectation function $\Phi$. Here, the bias of the rank estimators remains high even for small lags $\ell$ and large sample sizes $n$. This is due to the large differences between consecutive extreme values in the inputs $X_i$, amplified by the behavior of $\Phi : x \mapsto  x^2 - 3x$ (which is neither bounded nor Lipschitz in this case), causing the bias to remain high as $n \to \infty$. This example highlights the importance of the regularity conditions for the rank-based method to work and its potentially high bias, even when dealing with a single input.

\section{Conclusion}

The rank-based method proposed in \cite{GGKL20} provides an easily implementable estimator for first order Sobol indices. Specifically, given real-valued output $Y$ and input $X$, the second conditional moment $\eta = \mathbb E \big( \mathbb E(Y | X)^2 \big)$ is estimated by the lag-one cross-product of the outputs $Y_i$ ordered by increasing values of the inputs : 
$$ \widehat \eta^{(1)} = \frac 1 {n-1} \sum_{i=1}^{n-1} Y_{(i)} Y_{(i+1)} .    $$
Under regularity conditions on the expectation and variance of the response conditionally to the input, the estimator is known to be consistent and asymptotically Gaussian. In this paper, we discuss a natural extension of the method which consists in considering rank estimators obtained from higher order lags $\ell \geq 1$ :
$$ \widehat \eta^{(\ell)} = \frac 1 {n-\ell} \sum_{i=1}^{n-\ell} Y_{(i)} Y_{(i+\ell)}  \, , \, \ell =1,...,k.    $$
We show that these estimators share the same asymptotic properties under technical regularity conditions, provided that the maximal lag $k$ grows sufficiently slowly relative to $n$. We derive a closed form expression for the asymptotic covariance matrix of the collection $(\widehat \eta^{(1)}, ..., \widehat \eta^{(k)} )$, which allows to study the asymptotic behavior of the average estimator
$$ \widehat \eta^{(k)}_{\operatorname{av}} = \sum_{\ell = 1}^k \lambda_\ell \, \widehat \eta^{(\ell)}, $$
for suitable weights $\lambda_\ell, \ell = 1,...,k$. Base on the symmetry of the covariance matrix, the averaging procedure of \cite{lavancier2016general} justifies the equal weights $\lambda_\ell = 1/k$ as an asymptotically optimal choice. This is confirmed theoretically with the variance of average estimator $ \widehat \eta^{(k)}_{\operatorname{av}}$ reaching the efficiency bound of \cite{da2008efficient} for a regular estimator of $\eta$, whenever $k$ grows to infinity sufficiently slowly. In practice, the rule of thumb $k = \lfloor n^{1/3} \rfloor$ provides an entirely satisfactory choice in the various simulated examples, while verifying all the technical conditions for asymptotic efficiency. The theoretical results, as well as the importance of the regularity assumptions, are well validated by the numerical analysis. 

\section{Appendix} 
Let $\mathcal F_n$ denote the $\sigma$-algebra generated by $X_1,...,X_n$. The proofs of the results rely essentially on firstly investigating the distribution of the various estimators conditionally to $\mathcal F_n$. In particular, we exploit the fact that the $Y_{(i)}$'s remain independent conditionally to $X_1,...,X_n$ despite the sample re-shuffling, since the permutation that orders the inputs increasingly is $\mathcal F_n$-measurable.\\  
%To ease notation, we consider the sample $(Y_i, X_i), i=1,...,n$ as already ordered by increasing values of the $X_i$'s, keeping in mind that the sample is no longer iid in this case. 

To ease notation, we shall write $\phi_i = \Phi(X_i) = \mathbb E(Y_i | X_i) $ and $v_i = V(X_i) = \var(Y_i | X_i)$ for all $i=1,...,n$, and similarly for the ordered sample, e.g.~$\phi_{(i)}= \Phi(X_{(i)})$, $v_{(i)} = V(X_{(i)})$. 

\subsection*{Technical lemmas} 

\begin{lemma}\label{lem:espX} If \eqref{H1} and \eqref{H2} hold, then for $\ell = 1,...,k$, 
$$ \Big|  \mathbb E \big(\widehat \eta^{(\ell)} \, | \, \mathcal F_n \big) - \frac 1 n \sum_{i=1}^n \phi_i^2 \Big| \leq \frac \ell{n - \ell} \Big( L_\Phi M_\Phi \Delta_n + 2 M_\Phi^2  \Big). $$
%O \Big(\Delta_{k,n} + \frac{k}{n-k} \Big) .    $$
\end{lemma}

\textit{Proof.} Remark that  $ \mathbb E \big(Y_{(i)} Y_{(i+\ell)} \,  | \, \mathcal F_n \big) =  \phi_{(i)} \phi_{(i+\ell)}$ due to $Y_{(i)}$ and $Y_{(i+\ell)}$ being independent conditionally to $X_1,...,X_n$. It follows
%$$ \Big| \mathbb E \big(Y_i Y_{i+\ell} | \mathcal F_n \big) - \Phi^2(X_i)  \Big| = \Big|  \Phi(X_i) \big( \Phi(X_{i+\ell}) - \Phi(X_i) \big) \Big| \leq L_\Phi M_\Phi \Delta_{k,n} . $$
$$ \Big| \mathbb E \big(Y_{(i)} Y_{(i+\ell)} \, | \, \mathcal F_n \big) - \phi_{(i)}^2  \Big| = \Big|  \phi_{(i)} \big( \phi_{(i+\ell)} - \phi_{(i)} \big) \Big| \leq L_\Phi M_\Phi \big( X_{(i+\ell)} - X_{(i)} \big), $$
leading to
$$ \bigg| \mathbb E \big( \widehat \eta^{(\ell)} \, | \, \mathcal F_n \big)  - \frac 1 {n-\ell} \sum_{i=1}^{n-\ell} \phi_{(i)}^2 \bigg| 
%= \frac 1 {n-\ell}  \bigg|\sum_{i=1}^{n-\ell} \Big( \mathbb E \big( Y_i Y_{i+\ell} | \mathcal F_n \big)  - \Phi^2(X_i) \Big) \bigg| 
\leq L_\Phi M_\Phi  \frac {1} {n-\ell} \sum_{i=1}^{n-\ell} \big( X_{(i+\ell)} - X_{(i)} \big)  \leq L_\Phi M_\Phi \frac{\ell}{n-\ell} \Delta_n ,  $$
using Equation~\eqref{eq:teles}. Finally, summing over $n-\ell$ terms instead on $n$ deviates of at most
$$ \bigg| \frac 1 {n-\ell} \sum_{i=1}^{n-\ell} \phi_{(i)}^2 - \frac 1 n \sum_{i=1}^{n} \phi_i^2 \bigg| \leq \frac{2 \ell}{n-\ell} M_\Phi^2 $$ 
%= O \Big( \frac{k}{n-k} \Big)  $$
by \eqref{H1}, ending the proof. $\hfill \square$ \\

\begin{lemma}\label{lem:varX} Under Assumptions \eqref{H1} and \eqref{H2}, we have for $\ell < n$, 
% $$  \var \big( \widehat \eta^{(\ell)} | \mathcal F_n \big)  = \frac{1}{(n - \ell)} \frac 1 n \sum_{i=1}^n \Big( g_2^2(X_i) + 2 g_2(X_i) \Phi^2(X_i) - 3 \Phi^4(X_i) \Big) + o_P \Big( \frac{ 1} {n} \Big)    $$
$$ \bigg| (n - \ell)  \var \big( \widehat \eta^{(\ell)} \, | \, \mathcal F_n \big)  - \frac 1 n \sum_{i=1}^n \big( 4 \phi^2_i v_i + v^2_i \big) \bigg| \leq \frac{\ell}{n - \ell} \big( C_1 \Delta_n + C_2 \big) , $$
where $C_1, C_2$ are positive constants that depend only on $\Phi $ and $V$. \\
\end{lemma}

\textit{Proof.} 
%First remark that 
%$$ \mathbb E \big(\widehat \eta^{(\ell)2} \, | \, \mathcal F_n \big)  = \frac 1 {(n- \ell)^2} \sum_{i,j =1}^{n-\ell}  \mathbb E  \big( Y_{(i)} Y_{(i+\ell)} Y_{(j)} Y_{(j + \ell)} \, | \, \mathcal F_n \big) $$
%while 
%$$ \mathbb E \big( \widehat \eta^{(\ell)} \, | \, \mathcal F_n \big)^2 = \frac 1 {(n- \ell)^2}  \sum_{i, j=1}^{n - \ell}  \phi_{(i)} \phi_{(i+\ell)} \phi_{(j)} \phi_{(j + \ell)} . $$
Let $Z_{i,j,\ell} = \cov \big( Y_{(i)} Y_{(i+\ell)}, Y_{(j)} Y_{(j + \ell)} \, | \, \mathcal F_n \big)$
%$ = \mathbb E  \big( Y_{(i)} Y_{(i+\ell)} Y_{(j)} Y_{(j + \ell)} \, | \, \mathcal F_n \big) -  \phi_{(i)} \phi_{(i+\ell)} \phi_{(j)} \phi_{(j + \ell)} $$ 
and for a given $i \in \{ 1,...,n-\ell \}$, consider the set $S_{i,\ell} \subset \{1,...,n  \}$ of indices $j$ such that $i, i+ \ell, j,  j+\ell$ are not all distinct :
$$ S_{i,\ell}  = \{ i, i- \ell, i + \ell \} \cap \{ 1,...,n \} .  $$
Since $ Z_{i,j,\ell} = 0 $ if $ j \notin S_{i, \ell}$ by independence conditionally to $\mathcal F_n$, we have
$$ \var \big(\widehat \eta^{(\ell)} \, | \, \mathcal F_n \big) 
%= \mathbb E \big(\widehat \eta^{(\ell)2} \, | \, \mathcal F_n \big)  - \mathbb E \big( \widehat \eta^{(\ell)} \, | \, \mathcal F_n \big)^2  
%= \frac 1 {(n- \ell)^2} \sum_{i,j =1}^{n-\ell}  \mathbb E  \big( Y_{(i)} Y_{(i+\ell)} Y_{(j)} Y_{(j + \ell)} \, | \, \mathcal F_n \big) 
= \frac 1 {(n- \ell)^2} \sum_{i=1}^{n-\ell}  \sum_{j \in S_{i, \ell}} Z_{i,j,\ell} . $$
For $i=\ell+1,...,n-2\ell$, the set $S_{i, \ell}$ contains exactly three elements that can be dealt with separately:
\begin{itemize}
\item If $j = i$, $ Z_{i,j,\ell} = \phi_{(i)}^2 v_{(i + \ell)}  + \phi_{(i + \ell)}^2 v_{(i)} + v_{(i)} v_{(i + \ell)}$, and by \eqref{H1} and \eqref{H2}, 
$$  \Big| Z_{i,j,\ell}  - v_{(i)} \big( 2 \phi_{(i)}^2 + v_{(i)} \big) \Big| \leq 
(M_\Phi L_V + M_V L_\Phi + M_V L_V)   \big( X_{(i+\ell)} - X_{(i)} \big) . $$
\item If $j = i- \ell$ (and $i > \ell$), $ Z_{i,j,\ell} =  \phi_{(i-\ell)} \phi_{(i+ \ell)} v_{i} $ and
$$  \Big| Z_{i,j,\ell}  - \phi_{(i)}^2 v_{(i)} \Big| \leq M_\Phi M_V L_\Phi \big( X_{(i+\ell)} - X_{(i- \ell)} \big) . $$
\item If $j = i+ \ell$ (and $i \leq n-2\ell$), $ Z_{i,j,\ell} =  \phi_{(i)} \phi_{(i+ 2 \ell)} v_{(i+ \ell)} $ and
$$  \Big| Z_{i,j,\ell}  - \phi_{(i)}^2 v_{(i)}  \Big| \leq M_\Phi (M_V L_\Phi + M_\Phi L_V) \big( X_{(i+2\ell)} - X_{(i)} \big) . $$
\end{itemize}
Gathering all three terms, we obtain for all $i=\ell+1,...,n-2\ell$,
$$ \bigg| \sum_{j \in S_{i, \ell}} Z_{i,j,\ell} - \big( 4 \phi_{(i)}^2 v_{(i)} + v_{(i)}^2 \big) \bigg| \leq \frac{C_1}3 \big( X_{(i+2\ell)} - X_{(i- \ell)} \big)  $$
for some $C_1> 0$. Moreover, for the $2 \ell$ terms corresponding to $i \leq \ell$ and $ n - 2 \ell < i \leq n- \ell$, we can use the crude bound
$$ \bigg| \sum_{j \in S_{i, \ell}} Z_{i,j,\ell} - \big( 4 \phi_{(i)}^2 v_{(i)} + v_{(i)}^2 \big) \bigg| \leq 2(4  M_\Phi^2 M_V + M_V^2) . $$
Using the telescoping argument of Equation~\eqref{eq:teles}, we deduce
$$  \bigg| \sum_{i=1}^{n-\ell} \sum_{j \in S_{i, \ell}} Z_{i,j,\ell} - \sum_{i=1}^{n-\ell} \big( 4 \phi_{(i)}^2 v_{(i)} + v_{(i)}^2 \big) \bigg| \leq  \big( C_1 \Delta_n + 4(4  M_\Phi^2 M_V + M_V^2) \big) \ell . $$
The missing terms for $i > n - \ell$ can be bounded similarly by
$$ \bigg| \sum_{i=1}^{n-\ell} \big( 4 \phi_{(i)}^2 v_{(i)} + v_{(i)}^2 \big) - \sum_{i=1}^{n} \big( 4 \phi_{(i)}^2 v_{(i)} + v_{(i)}^2 \big)  \bigg| \leq \big( 4 M_\phi^2 M_V + M_V^2 \big) \ell  $$
and the result follows easily from here, using the triangular inequality. 
%$$  \bigg| \sum_{i=\ell+1}^{n-2\ell} \sum_{j \in S_{i, \ell}} Z_{i,j,\ell} - \sum_{i=\ell+1}^{n-2\ell} \big( 4 \phi_{(i)}^2 v_{(i)} + v_{(i)}^2 \big) \bigg| \leq C_1 \ell \Delta_n  $$
%Finally, because $Z_{i,j,\ell}$ is uniformly bounded by some positive constant $C_2$ in view of $\eqref{H1}$, the difference induced by summing over $i= \ell+1 $ to $n-2\ell$ instead of $i=1,...,n$ can be bounded by $C_2 \ell$ for some constant $C_2 > 0$, i.e.
%$$ \bigg|  \sum_{i=\ell+1}^{n-2\ell} \sum_{j \in S_{i, \ell}} Z_{i,j,\ell} -  \sum_{i=1}^{n- \ell} \sum_{j \in S_{i, \ell}} Z_{i,j,\ell} \bigg| \leq C_2 \ell  \big) $$
$ \hfill \square$

\subsection*{Proof of Proposition \ref{prop:var}} 
\noindent The result follows from Lemmas \ref{lem:espX} and \ref{lem:varX}, using the variance decomposition 
$$ \var \big( \widehat \eta^{(\ell)} \big)  = \mathbb E \Big( \var \big( \widehat \eta^{(\ell)} \, | \, \mathcal F_n  \big)  \Big) + \var \Big( \mathbb E \big( \widehat \eta^{(\ell)} | \mathcal F_n \big) \Big).  $$
%By assumption on $k$, $O(k/(n-k)) = o(1/\sqrt n)$ so that all such terms become negligible and can be safely ignored. 
Using Lemma \ref{lem:varX}, the first term in the variance decomposition is easily shown to verify 
$$ (n-\ell) \, \mathbb E \Big( \var \big( \widehat \eta^{(\ell)} | \mathcal F_n  \big)  \Big) =  4 \, \mathbb E \big( \Phi^2(X) \sigma^2(X) \big) + \mathbb E \big(\sigma^4(X) \big) + o (1) .$$
using that $\mathbb E(\ell \Delta_{n}) \leq \mathbb E(k \Delta_{n}) \leq \sqrt{\mathbb E(k^2 \Delta_{n}^2)} = o \big(1/\sqrt n \big) = o(1)$ by \eqref{H3}. For the second term, we have
$$ \var \Big( \frac 1 n \sum_{i=1}^n \phi_i^2 \Big) = \frac 1 n  \var \big( \Phi^2(X) \big)$$ %\leq \frac{M_\Phi^4}{n} $$
so that Lemma \ref{lem:espX} combined with \eqref{H3} give us directly
$$ \var \Big( \mathbb E \big( \widehat \eta^{(\ell)} | \mathcal F_n \big) \Big) = \frac{ \var \big( \Phi^2(X) \big)} n + o \Big( \frac 1 n \Big).   $$
%Let $B_{\ell,n} =  \mathbb E \big( \widehat \eta^{(\ell)} | \mathcal F_n \big) - \frac 1 n \sum_{i=1}^n \Phi^2 (X_i) $, we have by Cauchy-Schwarz's inequality, 
%$$ \Big| \var \Big( \mathbb E \big( \widehat \eta^{(\ell)} | \mathcal F_n \big) \Big) - \frac{ \var \big( \Phi^2(X) \big)} n \Big| \leq 2 \frac{M_\Phi^2 }{\sqrt n} \sqrt{\mathbb E(B_{\ell,n}^2)} + \mathbb E(B_{\ell,n}^2). $$
%We know from assumptions and the proof of Lemma \ref{lem:espX} that $ | B_{k,n} | \leq L_\Phi M_\Phi  \Delta_{k,n} + o\big(1/\sqrt n \big) $, leading to $\mathbb E \big(B_{\ell,n}^2 \big) = o \big(1/n \big)$. 
Hence,
$$ (n- \ell ) \var \Big( \mathbb E \big( \widehat \eta^{(\ell)} | \mathcal F_n \big) \Big) =  \frac{n-\ell} n  \var \big( \Phi^2(X) \big) + o (1) = \var \big( \Phi^2(X) \big) + o (1)   $$
and the result follows. $ \hfill \square$

\subsection*{Proof of Proposition \ref{prop:cov}}
\noindent We follow the same steps as in the proofs of Lemma \ref{lem:varX} and Proposition \ref{prop:var}, starting with 
$$ \cov \big(\widehat \eta^{(\ell)} , \widehat \eta^{(m)} \, | \, \mathcal F_n \big)  =  \frac 1 {(n- \ell)(n- m)} \sum_{i=1}^{n-\ell}  \sum_{j \in S_{i, \ell, m}} Z_{i,j,\ell, m} $$
where $Z_{i,j,\ell,m} = \mathbb E  \big( Y_{(i)} Y_{(i+\ell)} Y_{(j)} Y_{(j + m)} \, | \, \mathcal F_n \big) - \phi_{(i)} \phi_{(i+\ell)} \phi_{(j)} \phi_{(j + m)}$ and for all $i=1,..., n -\ell$, 
$$ S_{i,\ell,m} = \{ i, i+ \ell, i - m , i + \ell - m \} \cap \{1,...,n \} . $$
In the (at most) four cases for $j \in S_{i,\ell,m}$ and $i \in \{ m+1,..., n - \ell - m \}$, 
$$ \Big| Z_{i,j,\ell,m}  - \phi_{(i)}^2 v_{(i)} \Big| \leq C \big( X_{(i+\ell+m)} - X_{(i-m)} \big) $$
for some constant $C >0$. Using the same arguments, we arrive at
$$ \bigg| (n-\ell) \cov \big( \widehat \eta^{(\ell)}, \widehat \eta^{(m)} \, | \, \mathcal F_n \big) - \frac 4 n \sum_{i=1}^n \phi_i^2  v_i \bigg| \leq \frac{\ell}{n-\ell} \big( C_1'  \Delta_n + C_2' \big),$$
for some constants $C_1', C_2' > 0$. On the other hand
$$  \cov \Big( \mathbb E \big( \widehat \eta^{(\ell)} \, | \, \mathcal F_n \big) , \mathbb E \big( \widehat \eta^{(m)} \, | \, \mathcal F_n \big) \Big) =  \frac{ \var \big( \Phi^2(X) \big)} n + o \Big( \frac 1 n \Big) $$
by Lemma \ref{lem:espX}. We conclude using the decomposition
$$   \cov \big( \widehat \eta^{(\ell)}, \widehat \eta^{(m)} \big) = \mathbb E \Big( \cov \big(\widehat \eta^{(\ell)} , \widehat \eta^{(m)} \, | \, \mathcal F_n \big)  \Big) + \cov \Big( \mathbb E \big( \widehat \eta^{(\ell)} \, | \, \mathcal F_n \big) , \mathbb E \big( \widehat \eta^{(m)} \, | \, \mathcal F_n \big) \Big). \vspace*{-0.72cm} $$
$\hfill \square$ \vspace*{0.4cm}

\subsection*{Proof of Theorem \ref{th:av}}
\noindent Using Equation~\eqref{eq:sigma_asympt}, we obtain by straightforward calculation
\begin{equation}\label{eq:asympt_k} n \, \var \big( \widehat \eta^{(k)}_{\operatorname{av}} \big) = \frac {1} {k^2} \sum_{\ell, m = 1}^k n \, \cov \big( \widehat \eta^{(\ell)}, \widehat \eta^{(m)} \big) = \sigma^2_{\operatorname{opt}} + \frac 1 k  \mathbb E \big(V^2(X) \big) + o(1),
\end{equation}
after verification in the proofs that he residual terms $o(1)$ of Propositions \ref{prop:var} and \ref{prop:cov} are negligible uniformly for all $\ell, m \leq k$. The asymptotic efficiency follows for $k$ growing to infinity. $\hfill \square$

\bibliographystyle{abbrv}
\bibliography{ref}

\end {document}